\documentclass[12pt]{article}

\usepackage{amsmath,amssymb,amsthm,amscd}
\usepackage[colorlinks, linkcolor=blue, anchorcolor=gree, citecolor=red]{hyperref}
\usepackage[numbers,sort&compress]{natbib}
\usepackage{graphicx}

\begin{document}

\newcommand{\ts}{\,}
\newcommand{\tss}{\hspace{1pt}}
\newcommand{\Mat}{{\rm{Mat}}}
\newcommand{\CC}{\mathbb{C}}
\newcommand{\Sym}{\mathfrak S}

\newtheorem{thm}{Theorem}[section]
\newtheorem{pro}[thm]{Proposition}
\newtheorem{lem}[thm]{Lemma}
\newtheorem{cor}[thm]{Corollary}
\theoremstyle{definition}
\newtheorem{ex}[thm]{Example}
\newtheorem{remark}[thm]{Remark}
\newcommand{\bth}{\begin{thm}}
\renewcommand{\eth}{\end{thm}}
\newcommand{\bex}{\begin{examp}}
\newcommand{\eex}{\end{examp}}
\newcommand{\bre}{\begin{remark}}
\newcommand{\ere}{\end{remark}}

\newcommand{\bal}{\begin{aligned}}
\newcommand{\eal}{\end{aligned}}
\newcommand{\beq}{\begin{equation}}
\newcommand{\eeq}{\end{equation}}
\newcommand{\ben}{\begin{equation*}}
\newcommand{\een}{\end{equation*}}

\newcommand{\bpf}{\begin{proof}}
\newcommand{\epf}{\end{proof}}
\renewcommand{\thefootnote}{}

\def\beql#1{\begin{equation}\label{#1}}
\title{\Large\bf  Power graphs of (non)orientable genus two}

\author{{\sc Xuanlong Ma$^1$, Gary L. Walls$^2$, Kaishun Wang$^1$
}\\[15pt]
{\small $^1$Sch. Math. Sci. {\rm \&} Lab. Math. Com. Sys., Beijing Normal University, Beijing, 100875, China}\\
{\small $^2$Department of Mathematics, Southeastern Louisiana University, Hammond, LA 70402, USA}\\
}

 \date{}

\maketitle

\begin{abstract}
The power graph $\Gamma_G$ of a finite group $G$ is the graph whose vertex set is the group, two distinct elements being adjacent if one is a power of the other. In this paper, we  classify the finite groups whose power graphs have (non)orientable genus two.
\end{abstract}



\footnote{E-mail addresses: xuanlma@mail.bnu.edu.cn (X. Ma), gary.walls@selu.edu (G.L. Walls) wangks@bnu.edu.cn (K. Wang).}
\section{Introduction}
Throughout this paper, every graph is finite, simple and connected.
A graph $\Gamma$ is called a {\em planar graph} if $\Gamma$ can be drawn in the plane so that no two of its edges cross each other. In addition, we say that $\Gamma$ can be embedded in the plane.
A non-planar graph can be  embedded in some surface obtained from the sphere by attaching some handles or crosscaps.
We denote by $\mathbb{S}_k$ a sphere with $k$ handles and by $\mathbb{N}_k$ a sphere with $k$ crosscaps.
Note that both $\mathbb{S}_0$ and $\mathbb{N}_0$ are the sphere itself, and $\mathbb{S}_1$
and $\mathbb{N}_1$ are a torus and a projective plan, respectively.
The smallest non-negative integer $k$
such that a graph $\Gamma$ can be embedded on $\mathbb{S}_k$ is called the {\em orientable genus} or {\em genus} of $\Gamma$, and is denoted by $\gamma(\Gamma)$. The {\em nonorientable genus} of $\Gamma$, denoted by $\overline{\gamma}(\Gamma)$, is the smallest integer $k$ such that $\Gamma$ can be embedded on $\mathbb{N}_k$.

The problem of finding the graph genus is NP-hard \cite{T89}.
Many research articles have appeared on the genus of graphs constructed from some algebraic structures. For example, Wang \cite{W06} found all rings of two specific forms that have genus at most one. Chiang-Hsieh et al. \cite{CSW10} characterized the commutative rings of genus one.
Very recently, Rajkumar and Devi \cite{RD15} classified the finite groups whose intersection graphs of subgroups have (non)orientable genus one.
Afkhami et al. \cite{A15} classified planar, toroidal, and projective commuting and noncommuting graphs of finite groups.

Here we study the genus of the power graph of a finite group.
The {\em undirected power graph} $\Gamma_G$ of a group $G$ has the vertex set $G$ and two distinct elements are adjacent if one is a power of the other. The concept of a power graph was first introduced and considered by Kelarev and Quinn \cite{KQ00}.
 Note that since the paper deals only with undirected graphs, for convenience,  throughout we use the term ``power graph" to refer to an undirected power graph defined as above.
Recently, many interesting results on power graphs have been obtained, see \cite{Cam,CGh,CGS,Kel2,Kel21,Kel22,FMW,FMW1,Ma,MFW}. Furthermore, \cite{AKC13} is a survey that is a detailed list of results and open questions on power graphs.

Doostabadi and Farrokhi D.G. \cite{D14} classified the finite groups whose power graphs have (non)orientable genus one.
The goal of the paper is to find all power graphs of (non)orientable genus two.
Our main results are the following theorems.

\begin{thm}\label{mainthm1}
Let $G$ be a finite group. Then $\Gamma_{G}$ has orientable genus two
if and only if $G$ is isomorphic to one group in Table \ref{tab1}.
\begin{table}[htbp]\centering
\caption{\small All finite groups $G$  with  $\gamma(\Gamma_G)=2$\label{tab1}}
\vskip2mm
{\footnotesize\tabcolsep=10pt
\begin{tabular}{cc}
\hline
GAP ID    &  Group                                  \\  \hline
$[8,1]$   & $\mathbb{Z}_8$                           \\
$[12,5]$  & $\mathbb{Z}_2\times \mathbb{Z}_6$       \\
$[16,7]$  & $D_{16}$                                \\
$[16,8]$  & $QD_{16}$                               \\
$[16,9]$  & $Q_{16}$                                 \\
$[18,3]$  & $\mathbb{Z}_3\times S_3$                 \\
$[24,7]$  & $\mathbb{Z}_2\times (\mathbb{Z}_3\rtimes_{\varphi}\mathbb{Z}_4)$  \\
$[24,8]$  & $(\mathbb{Z}_6\times\mathbb{Z}_2)\rtimes_{\varphi}\mathbb{Z}_2$   \\
$[24,14]$ & $\mathbb{Z}_2\times\mathbb{Z}_2 \times S_3$                        \\
$[36,11]$  & $\mathbb{Z}_3\times A_4$                                          \\
$[72,43]$  & $(\mathbb{Z}_3\times A_4)\rtimes_{\varphi}\mathbb{Z}_2$           \\
 \hline
\end{tabular}}
\end{table}
\end{thm}

\begin{thm}\label{mainthm2}
There is no finite group $G$ such that $\Gamma_G$ has nonorientable genus two.
\end{thm}

\section{Preliminaries}\label{sec:}
In this section we briefly recall some notation, terminology, and basic results and prove
a lemma  which we need in the sequel.

Let $\Gamma$ be a graph.
Denote by $V(\Gamma)$ and $E(\Gamma)$ the vertex set and the edge set of $\Gamma$, respectively.
An edge of $\Gamma$ is denoted simply by $ab$
where $a,b\in V(\Gamma)$.
If $V'\subseteq V(\Gamma)$, we define $\Gamma - V'$ to be the subgraph of $\Gamma$ obtained by deleting the vertices in $V'$ and all edges incident with them.
Similarly,
if  $E'\subseteq E(\Gamma)$, then $\Gamma - E'$ is the subgraph of $\Gamma$ obtained by deleting
the edges in $E'$.
For two vertex-disjoint graphs $\Gamma$ and $\Delta$, $\Gamma \cup \Delta$ denotes the
graph with vertex set $V(\Gamma)\cup V(\Delta)$ and edge set $E(\Gamma)\cup E(\Delta)$, and
$\Gamma + \Delta$ consists of $\Gamma \cup \Delta$ and all edges joining a vertex of
$\Gamma$ and a vertex of $\Delta$. A union of $k$ isomorphic graphs $\Gamma$ is denoted by $k\Gamma$.
We use the natation $\lceil x\rceil$ to denote the least integer that is greater than or
equal to $x$.
Denote by $K_n$ and $K_{m,n}$ the complete graph of order $n$ and the complete bipartite graph,
respectively.

For any subgraph $\Delta$ of a graph $\Gamma$, one easy observation is that
$\gamma(\Delta)\le \gamma(\Gamma)$ and $\overline{\gamma}(\Delta)\le \overline{\gamma}(\Gamma)$.
The following result gives the (non)orientable genus of a complete graph and a complete bipartite graph.

\begin{thm}{\rm (\cite[p. 58, p. 152]{Whi})}\label{ccgenus}
Let $n$ be an integer at least $3$. Then

$(a)$ $\gamma(K_n)=\lceil \frac{1}{12}(n-3)(n-4)\rceil$.

$(b)$  $\overline{\gamma}(K_n)=\lceil \frac{1}{6}(n-3)(n-4)\rceil$ if $n\ne 7$;  $\overline{\gamma}(K_7)=3$.

$(c)$ $\gamma(K_{m,n})=\lceil \frac{1}{4}(m-2)(n-2)\rceil$.

$(d)$ $\overline{\gamma}(K_{m,n})=\lceil \frac{1}{2}(m-2)(n-2)\rceil$.
\end{thm}



A {\em block} of a graph $\Gamma$ is a maximal connected subgraph $B$ of $\Gamma$  with respect to the property that removing any vertex of $B$ does not disconnect $B$.
The following result tells us how to compute the (non)orientable genus of a graph by its blocks.

\begin{thm}{\rm (\cite[Theorem 1]{BHKY62}, \cite[Corollary 3]{SBblock})}\label{moo}
Let $\Gamma$ be a connected graph with $n$ blocks $B_1, \cdots, B_n$. Then

$(1)$ $$\gamma(\Gamma)=\sum_{i=1}^n\gamma(B_i).$$

$(2)$ If $\overline{\gamma}(B_i)= 2\gamma(B_i)+1$ for each $i$, then $$\overline{\gamma}(\Gamma)=1-n+\sum_{i=1}^n\overline{\gamma}(B_i).$$
Otherwise,
$$\overline{\gamma}(\Gamma)=2n-\sum_{i=1}^n\mu(B_i),$$
where $\mu(B_i)=\max\{2-2\gamma(B_i), 2-\overline{\gamma}(B_i)\}$.
\end{thm}

Now we state Euler's formula for $\mathbb{S}_k$ and $\mathbb{N}_k$.
\begin{thm}\label{ef}
Suppose that $\Gamma$ is a connected graph that is embedded on a surface $S$, resulting in $f$ faces. Then

$(a)$ If $S=\mathbb{S}_{\gamma(\Gamma)}$, then $|V(\Gamma)|-|E(\Gamma)|+f=2-2\gamma(\Gamma)$.

$(b)$ If $S=\mathbb{N}_{\overline{\gamma}(\Gamma)}$, then $|V(\Gamma)|-|E(\Gamma)|+f=2-\overline{\gamma}(\Gamma)$.
\end{thm}

All groups considered in this paper are finite. An element of order $2$ in a
group is called an {\em involution}.
Let $G$ be a group and $g$ be an element of $G$.
Denote by $|G|$ and $|g|$ the orders of $G$ and  $g$, respectively.
Let $S$ be a subset of $G$.
The set of prime divisors of $|S|$ is denoted by $\pi(S)$ and
the set of natural numbers consisting of orders of elements of $S$ is denoted by $\pi_e(S)$.
Also $\mathbb{Z}_n^m$ is used for the $m$-fold direct product of the cyclic group $\mathbb{Z}_n$ with itself.
Denote by $\mathbb{Z}_n$ and  $D_{2n}$ the cyclic group of order $n$ and the dihedral group of order $2n$, respectively.

Planar power graphs were characterized in \cite{MAN}.

\begin{thm}{\rm (\cite[Corollary 4]{MAN})}\label{plan}
Let $G$ be a group. Then $\Gamma_G$ is planar if and only if $\pi_e(G)=\{1,2,3,4\}$.
\end{thm}

For a subset $S$ of a group $G$, $\Gamma_{G[S]}$ denotes the induced subgraph of $\Gamma_G$ by $S$ and if the situation is unambiguous, then we denote $\Gamma_{G[S]}$ simply by $\Gamma_{S}$.

\begin{lem}\label{condi}
Let $G$ be a group. Suppose that $S$ is a union of some subgroups of $G$ such that $\pi_e(G\setminus S)\subseteq\{2,3,4\}$. If $\gamma(\Gamma_S)=r$ {\rm (}resp.  $\overline{\gamma}(\Gamma_S)=r${\rm)}, then $\gamma(\Gamma_G)=r$ {\rm (}resp.  $\overline{\gamma}(\Gamma_G)=r${\rm)}.
\end{lem}
\bpf
We first prove that if $\gamma(\Gamma_S)=r$, then $\gamma(\Gamma_G)=r$.
Since $\gamma(\Gamma_S)=r$, $\Gamma_S$ can be embedded  on $\mathbb{S}_r$. Now
we fix an embedding $\mathbb{E}$ of  $\Gamma_S$ on $\mathbb{S}_r$. Clearly, if
$G\setminus S=\emptyset$, then $\gamma(\Gamma_G)=r$. Thus, we may assume that $G\setminus S\ne\emptyset$.

\medskip
\noindent {\bf Case 1.} $\Gamma_G$ has no edge $ab$ such that $a\in S\setminus\{e\}$ and
$b\in G\setminus S$, where $e$ is the identity element of $G$.
\medskip

Take a face $F$ containing
$e$ in $\mathbb{E}$. By Theorem \ref{plan}, we may insert all vertices in $G\setminus S$ and all edges incident with them in $F$ without crossings.
This implies that $\Gamma_G$ also can embed on a surface of genus $r$ and so $\gamma(\Gamma_G)=r$, as required.

\medskip
\noindent {\bf Case 2.}  $\Gamma_G$ has an edge $ab$ such that $a\in S\setminus\{e\}$ and $b\in G\setminus S$.
\medskip

It is easy to
see that $a$ and $b$ are an involution and an element of order $4$, respectively.
Let $b,b^{-1},b_1,b_1^{-1},\cdots,b_t,b_t^{-1}$ be all elements that are adjacent to
$a$ in $G\setminus S$. Note that $|b_i|=4$ for each $i=1,\ldots,t$.
Since $e$ and $a$ are adjacent, we may
take a face $F$ containing edge $ea$ in $\mathbb{E}$.
Note that for any vertex of $\{b,b^{-1},b_1,b_1^{-1},\cdots,b_t,b_t^{-1}\}$, there are exactly
three vertices that are adjacent to it.
Then we can insert the vertices $b,b^{-1},b_1,b_1^{-1},\cdots,b_t,b_t^{-1}$ and and all edges incident with them in $F$ without crossings, as shown in Figure \ref{aad1}.
\begin{figure}[hptb]
  \centering
  \includegraphics[width=7cm]{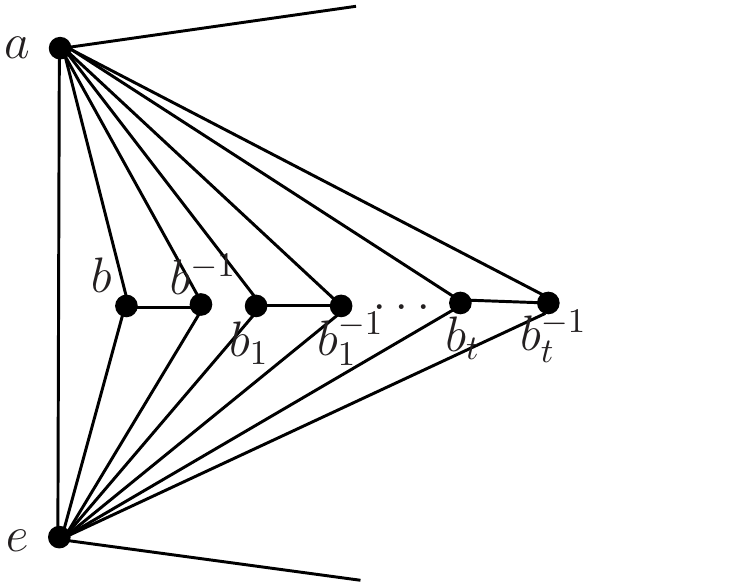}\\
  \caption{Insert some vertices and edges in face $F$.}\label{aad1}
\end{figure}
It follows that if $G$ has other involutions in $S\setminus\{e\}$ that are adjacent to some elements in $G\setminus S$, then we may insert them in some faces of $\mathbb{E}$ without crossings.
For the remainder vertices and edges in $\Gamma_{G\setminus S}$, by Theorem \ref{plan} we also can insert them in a face containing $e$ of $\mathbb{E}$ without crossings.
This implies that $\gamma(\Gamma_G)=r$, as desired.

Similarly, we have that if $\overline{\gamma}(\Gamma_S)=r$, then $\overline{\gamma}(\Gamma_G)=r$.
\epf

\section{Groups with some cyclic subgroups of order $6$}\label{Groups6}

In this section we prove some results on finite groups with some cyclic subgroups of order $6$.

\begin{lem}\label{two6}
There is no group that has precisely two cyclic subgroups of order $6$.
\end{lem}
\bpf
Suppose, towards a contradiction, that
there exists a group $G$ that has exactly two cyclic subgroups of order $6$, say $\langle a\rangle$ and $\langle b\rangle$. Note that $G$ has precisely four elements $a,a^{-1},b,b^{-1}$ of order $6$. Then $b^a=a,a^5,b$ or $b^5$, where $b^a$ is the conjugate of $b$ by $a$, that is, $b^a=a^{-1}ba$.
Since $\langle a\rangle\ne \langle b\rangle$, one has that $b^a= b$ or $b^{-1}$. It follows that $\langle a\rangle\langle b\rangle$ is a subgroup of order $36$, $18$ or $12$.
We check all groups of order $36$, $18$ and $12$ (for example using the
computer algebra system GAP \cite{G13}) and find that
there is no such group that has only four elements of order $6$, which is a contradiction.
\epf

\begin{lem}\label{claim}
Let $G$ be a group with $\pi_e(G)\subseteq\{1,2,3,4,6\}$. Suppose that
$G$ has precisely three cyclic subgroups. If there exist two cyclic subgroups of order $6$ of $G$ such that their intersection has order $3$, then the intersection of
any two cyclic subgroups of order $6$ of $G$ is of order $3$.
\end{lem}
\bpf
Let $\langle a\rangle$, $\langle b\rangle$, $\langle c\rangle$ be the three cyclic subgroups
of order $6$ of $G$ and assume without loss of generality that
$\langle b\rangle\cap \langle c\rangle=\langle b^2 \rangle =\langle c^2 \rangle$.

If $\langle a\rangle$ is a normal subgroup of $G$, then $[G: C_G(a)]\le 2$ by the $N/C$ Theorem (cf. \cite[Theorem 1.6.13]{Rob}, here $[G: C_G(a)]$ denotes the index of
$C_G(a)$ in $G$). It follows that $b^2\in C_G(a)$. Hence, $|ab^2|=6$.
If $ ab^2\in \langle a\rangle$, then $b^2\in \langle a\rangle$ so $b^2=a^2$ and we have
$a^2=b^2=c^2$, as required. If $ ab^2\in \langle b\rangle$, then $a\in \langle b\rangle$,
a contradiction. Similarly, if $ ab^2\in \langle c\rangle$, then again $ a\in \langle c\rangle$, since $b^2=c^2$ and this is a contradiction. The result follows in this case.

So assume that $\langle a\rangle$ is not normal in $G$. Then, without loss of generality
there is a $g\in G$, so that $\langle a\rangle^g=\langle b\rangle$. Now
$\langle a\rangle^g \cap \langle c\rangle=\langle b^2\rangle$. It follows that
$\langle a\rangle \cap \langle c\rangle^{g^{-1}}$ has order $3$.
Now $\langle c^{g^{-1}}\rangle$ cannot be $\langle a\rangle$, so it must be either $\langle b\rangle$ or $\langle c\rangle$. It follows that either $\langle a^2\rangle=\langle b^2\rangle$
or $\langle a^2\rangle=\langle c^2\rangle$ and result follows.
\epf

In GAP \cite{G13},
the GAP ID $[n,m]$ which
is a label that uniquely identifies the group in GAP, the first number $n$ in the square brackets is the order of the group, and the second number $m$ simply enumerates different groups of the same order.

\begin{thm}\label{gtheory3}
Let $G$ be a group with $\pi_e(G)\subseteq\{1,2,3,4,6\}$. Suppose that
$G$ has precisely three cyclic subgroups $\langle a\rangle$, $\langle b\rangle$, $\langle c\rangle$ of order $6$. Then there exist two cyclic subgroups of order $6$ in $G$ such that
their intersection has order $3$
if and only if $G$ is isomorphic to one group in Table \ref{tab2}.
\begin{table}[htbp]\centering
\caption{\small All finite groups satisfying the conditions \label{tab2}}
\vskip2mm
{\footnotesize\tabcolsep=10pt
\begin{tabular}{cccccc}
\hline
GAP ID    &  Group                              &    $\pi_e(G)$            \\  \hline
$[12,5]$  & $\mathbb{Z}_2\times \mathbb{Z}_6$   &    $\{1,2,3,6\}$         \\
$[18,3]$  & $\mathbb{Z}_3\times S_3$            &    $\{1,2,3,6\}$        \\
$[24,7]$  & $\mathbb{Z}_2\times (\mathbb{Z}_3\rtimes_{\varphi}\mathbb{Z}_4)$ &$\{1,2,3,4,6\}$ \\
$[24,8]$  & $(\mathbb{Z}_6\times\mathbb{Z}_2)\rtimes_{\varphi}\mathbb{Z}_2$ &$\{1,2,3,4,6\}$\\
$[24,14]$ & $\mathbb{Z}_2\times\mathbb{Z}_2 \times S_3$ &$\{1,2,3,6\}$ \\
$[36,11]$  & $\mathbb{Z}_3\times A_4$ &$\{1,2,3,6\}$\\
$[72,43]$  & $(\mathbb{Z}_3\times A_4)\rtimes_{\varphi}\mathbb{Z}_2$ &$\{1,2,3,4,6\}$\\
 \hline
\end{tabular}}
\end{table}
\end{thm}
\bpf
We check all groups of order at most $144$ (for example using the
computer algebra system GAP \cite{G13}) and find that the  groups $G$ in Table \ref{tab2} are precisely the groups of order at most $144$  satisfying the conditions:

$(a)$ $\pi_e(G)\subseteq\{1,2,3,4,6\}$.

$(b)$ $G$ has precisely three cyclic subgroups of order $6$.

$(c)$ There exist two cyclic subgroups of order $6$ of $G$ such that their intersection has order $3$.

In order to complete the proof,
next we prove if $G$ is a group satisfying the conditions, then $|G|\le 144$.

By Lemma \ref{claim}, we have that $\langle b\rangle\cap \langle c\rangle=\langle a\rangle\cap \langle b\rangle=\langle a\rangle\cap \langle c\rangle=\langle a^2\rangle=\langle b^2\rangle=\langle c^2\rangle$. Thus,
without loss of generality, we may assume that $a^2=b^2=c^2$.

\medskip
\noindent {\bf Case 1.} $a^3\in Z(G)$, the center of $G$.
\medskip

Note that for any  prime divisor $p$ of $|G|$, in general, the number of subgroups of order $p$ of $G$ is congruent to $1$ modulo $p$.
Since $a^3\in Z(G)$ and $G$ has exactly three cyclic subgroups  of order $6$, $G$ has a unique subgroup $\langle a^2\rangle$ of order $3$. This implies that $\langle a^2\rangle$ is normal in $G$.

Suppose that $a^2\in Z(G)$. Since $\pi_e(G)\subseteq\{1,2,3,4,6\}$, $G$ has no elements of order $4$. If $G$ has a unique involution, then it is easy to see that $G\cong \mathbb{Z}_6$, a contradiction. Therefore,
$G$ has exactly three involutions.
It follows that $G=\mathbb{Z}_2 \times \mathbb{Z}_6$.

Now suppose that $a^2\notin Z(G)$. Then $C_G(a^2)\ne G$. So,
by $N/C$ Theorem  we see that $|G|=2|C_G(a^2)|$.
Since $a^2=b^2=c^2$, we have $a,b,c\in C_G( a^2)$. Thus, $C_G( a^2)$ also satisfies the conditions $(a)$--$(c)$.  Moreover, note that $a^3, a^2\in Z(C_G(a^2))$ and the proof above, it follows that $C_G(a^2)\cong \mathbb{Z}_2 \times \mathbb{Z}_6$. Consequently, in this case we have $|G|=24$, as desired.

\medskip
\noindent {\bf Case 2.} $a^3\notin Z(G)$.
\medskip

By Case $1$ we may assume that $b^3,c^3\notin Z(G)$. Let $S=\{a^3,b^3,c^3\}$. It
follows that $S$ is a conjugacy class of $G$. Thus, we obtain that
$[G: C_G(a^3)]=[G: C_G(b^3)]=[G: C_G(c^3)]=3$. Furthermore, it is clear that $C_G(a^3)\cap C_G(b^3)\subseteq C_G(c^3)$.
Hence, one has that $C_G(a^3)\cap C_G(b^3)\cap C_G(c^3)=C_G(a^3)\cap C_G(b^3)$.
Similarly, one can get $C_G(a^3)\cap C_G(c^3)=C_G(a^3)\cap C_G(b^3)=C_G(b^3)\cap C_G(c^3)$.

Now let $\eta$ be the group action of $G$ on $S$, where each element of $G$ acts
on the conjugacy class $S$ by conjugation. Then
$${\rm Ker}(\eta)=C_G(a^3)\cap C_G(b^3)\cap C_G(c^3)=C_G(a^3)\cap C_G(b^3).$$
Note that $\eta$ is a homomorphism from $G$ to $S_3$. So
we have $[G: C_G(a^3)\cap C_G(b^3)]\le 6$.
It follows that $[C_G(a^3): C_G(a^3)\cap C_G(b^3)]\le 2$.

\medskip
\noindent {\bf Subcase 2.1.} $[C_G(a^3): C_G(a^3)\cap C_G(b^3)]=1$.
\medskip

Then
$[C_G(b^3): C_G(c^3)\cap C_G(b^3)]=1$, $[C_G(c^3): C_G(c^3)\cap C_G(a^3)]=1$
and $C_G(a^3)$ is normal in $G$. This implies that $C_G(a^3)=C_G(b^3)=C_G(c^3)$ and so $b^3,c^3\in C_G(a^3)$. Since $a^2=b^2=c^2$, one has that $a,b,c\in C_G(a^3)$.
This implies that
$C_G(a^3)$ is a group satisfying the conditions $(a)$--$(c)$. Note that $a^3\in Z(C_G(a^3))$. By Case $1$ we have $|C_G(a^3)|\le 24$. It follows that $|G|=3|C_G(a^3)|\le 72$, as desired.

\medskip
\noindent {\bf Subcase 2.2.} $[C_G(a^3): C_G(a^3)\cap C_G(b^3)]=2$.
\medskip

Then
$[C_G(b^3): C_G(b^3)\cap C_G(c^3)]=[C_G(c^3): C_G(c^3)\cap C_G(a^3)]=2$. Note that
$C_G(a^3)\cap C_G(b^3)$ is clearly  a normal subgroup of $G$.

Suppose that
$a^3\in C_G(a^3)\cap C_G(b^3)$. Then $S\subseteq C_G(a^3)\cap C_G(b^3)$ and $a^2\in C_G(a^3)\cap C_G(b^3)$. It follows that $a,b,c\in C_G(a^3)\cap C_G(b^3)$. So $C_G(a^3)\cap C_G(b^3)$ is
a group satisfying the conditions $(a)$--$(c)$ and $a^3\in Z(C_G(a^3)\cap C_G(b^3))$. In view of Case $1$,
we have that $|C_G(a^3)\cap C_G(b^3)|\le 24$. Thus, $|G|=6|C_G(a^3)\cap C_G(b^3)|\le 144$, as desired.

Suppose that $a^3\notin C_G(a^3)\cap C_G(b^3)$. Note that $S$ is a conjugacy class.
Then $c^3,b^3\notin C_G(a^3)\cap C_G(b^3)$. Note that $C_G(b^3)\ne C_G(a^3)$.
Thus, $a^3b^3\notin C_G(a^3)$. Now it is easy  to check that $|a^3b^3|=3$. This implies that
$\langle a^3,b^3\rangle\cong S_3$. Note that $\langle a^3,b^3\rangle \cap
C_G(a^3)\cap C_G(b^3)=1$ and $C_G(a^3)\cap C_G(b^3)$ is normal in $G$. We have
$G\cong \langle a^3,b^3\rangle \times (C_G(a^3)\cap C_G(b^3))$.
Since $G$ has precisely three cyclic subgroups of order $6$,
$C_G(a^3)\cap C_G(b^3)$ has no involutions and has only two elements of order $3$.
So $C_G(a^3)\cap C_G(b^3)\cong \mathbb{Z}_3$. This implies that
$G\cong S_3 \times \mathbb{Z}_3$, as desired.
\epf

\begin{thm}\label{gtheory4}
Let $G$ be a group with $\pi_e(G)\subseteq\{1,2,3,4,6\}$. Suppose that $G$ has precisely four cyclic subgroups of order $6$.
Then there is no finite group $G$ such that the intersection of two cyclic subgroups of order $6$ of $G$ has order $3$, and the intersection of the remaining two cyclic subgroups of order $6$ also has order $3$.
\end{thm}
\bpf
Suppose, for a contradiction, that $G$ is a such group of minimal order.
Let $\langle a\rangle, \langle b\rangle, \langle c\rangle, \langle d\rangle$ be the four cyclic subgroups of order $6$ of $G$, and
$|\langle a\rangle \cap \langle b\rangle|=3$ and $|\langle c\rangle\cap \langle d\rangle|=3$.
Without loss of generality, we may assume that $a^2=b^2$ and $c^2=d^2$. Then
$\langle a^3,b^3\rangle\subseteq C_G(a^2)$. Note that $\langle a^3,b^3\rangle$ has at least $3$
involutions. Thus, $\langle a,b\rangle$ has at least $3$ cyclic subgroups of order $6$.
So, without loss of generality, we may assume that $c\in \langle a,b\rangle$.

If  $d\in \langle a,b\rangle$, then it is clear that $\langle a,b,c,d\rangle \subseteq C_G(a^2)$.
Suppose that $d\notin \langle a,b\rangle$. Then $\langle a,b\rangle$ has precisely three cyclic subgroups of order $6$ and $|\langle a\rangle \cap \langle b\rangle|=3$.
By Lemma \ref{claim}, we may assume that $a^2=b^2=c^2$. Hence, we have $a^2=d^2$. It follows that $d\in C_G(a^2)$ and so
$\langle a,b,c,d\rangle \subseteq C_G(a^2)$. This implies that
we always may assume that $\langle a,b,c,d\rangle \subseteq C_G(a^2)$.

Now note that $G$ is of minimal order. We have $G= C_G(a^2)$. This implies that $a^2\in Z(G)$.
Since $G$ has precisely four cyclic groups of order $6$, $G$ has at most four involutions.
Now it is easy to check that the number of involutions of $G$ is $2$ or $4$.
This contradicts the fact that the number of involutions of a finite group of even order
is odd.
\epf

\section{Proof of Theorem \ref{mainthm1}}\label{sec:}

\begin{lem}\label{cycg} $\gamma(\Gamma_{\mathbb{Z}_n})=2$ if and only if $n=8$. In particular,
$\gamma(\Gamma_{\mathbb{Z}_n})\ge 3$ for $n\ge 9$.
\end{lem}
\bpf
Note that for a group $G$, $\Gamma_{G}$ is complete if and only if $G$ is a cyclic group of prime power order by \cite[Theorem 2.12]{CGS}. Thus if $n=8$, then $\gamma(\Gamma_{\mathbb{Z}_8})=2$ by Theorem \ref{ccgenus}.

Now suppose that $\gamma(\Gamma_{\mathbb{Z}_n})=2$. Since $\gamma(K_7)=1$, one has
$\gamma(\Gamma_{\mathbb{Z}_n})\le 1$ if $n\le 7$. Thus, we may assume that $n\ge 8$.
If $n\ge 9$, by \cite[Theorem 2]{MAN} one easy calculation shows that
the clique number of $\Gamma_{\mathbb{Z}_n}$ is greater than $8$ , and so its genus is not two by Theorem \ref{ccgenus}. This implies that $n=8$.
\epf

The generalized quaternion group $Q_{16}$ of order $16$ which is given by
$
\langle x,y: x^4=y^2, x^{8}=1, y^{-1}xy=x^{-1}\rangle.
$
$QD_{16}$ denotes the group of order $16$ that is given by $\langle a,b: a^8=b^2=1,bab=a^3\rangle$.

\begin{lem}\label{p2g}
Let $G$ be a $p$-group, where $p$ is a prime. Then
$\gamma(\Gamma_{G})=2$ if and only if $G$ is isomorphic to one of
the following groups:
\begin{equation}\label{2pg}
\mathbb{Z}_8,~D_{16},~Q_{16},~QD_{16}.
\end{equation}
\end{lem}
\bpf
By verifying we know that every group $G$ in (\ref{2pg}) has a unique cyclic subgroup $S$ of order $8$,
and $\pi_e(G\setminus S)\subseteq\{2,3,4\}$. Note that $\gamma(\Gamma_S)=\gamma(K_8)=2$.
By Lemma \ref{condi}, we have $\gamma(\Gamma_{G})=2$.

Now we suppose that $\gamma(\Gamma_{G})=2$.
By Lemma \ref{cycg}, every element of $G$ has order at most $8$. Thus, we
may assume that $p\le 7$.

Suppose $p=7$. Note that $\gamma(\Gamma_{\mathbb{Z}_7})= 1$. Then $|G|=7^n$ for some $n\ge 2$. This implies that
$G$ has a subgroup $A$ isomorphic to $\mathbb{Z}_7\times \mathbb{Z}_7$ or $\mathbb{Z}_{49}$.
Since $G$ has no elements of order $49$, $A\cong \mathbb{Z}_7\times \mathbb{Z}_7$.
So $\Gamma_A$ is isomorphic to $K_1+8K_6$ that has genus $8$ by Theorems \ref{moo}, a contradiction.
Similarly, we obtain $p\ne 5$.

If $p=3$, then $\pi_e(G)=\{1,3\}$, a contradiction by Theorem \ref{plan}. Thus,
now we may suppose that $|G|=2^n$ for some $n\ge 3$, and so $\pi_e(G)=\{1,2,4,8\}$.
If $G$ has two distinct cyclic subgroups of order $8$, then the subgraph of $\Gamma_G$ induced by the two cyclic subgroups contains a subgraph isomorphic to $K_1+(K_7\cup K_4)$ which has genus $3$ by Theorem \ref{moo}, a contradiction.
This implies that $G$ has a unique cyclic subgroup of order $8$, which is normal in $G$. Let $g\in G$ with $|g|=8$. If there exists an element $x$ in $G\setminus \langle g\rangle$ such that $x\in C_G( g)$, then $G$ has a subgroup isomorphic to $\mathbb{Z}_2\times \mathbb{Z}_8$, which is a contradiction since
$\mathbb{Z}_2\times \mathbb{Z}_8$ has two cyclic subgroups of order $8$.
Therefore, we have that $C_G(g)=\langle g\rangle$. This implies that the quotient group $G/\langle g\rangle$ is isomorphic to a subgroup of the full automorphism group of $\langle g\rangle$. Thus, we conclude that $|G|=8$, $16$ or $32$.
If $|G|=8$, then $G\cong \mathbb{Z}_8$. If $|G|=16$, then $G\cong D_{16}$, $Q_{16}$ or $QD_{16}$. Since there is no such  group $G$ of order $32$ such that $\pi_e(G)=\{1,2,4,8\}$ and $G$ has a unique cyclic subgroup of order $8$, we get the desired result.
\epf

\begin{lem}\label{three6}
Let $G$ be a group with $\pi_e(G)\subseteq\{1,2,3,4,6\}$. Suppose that $G$ has precisely three cyclic subgroups of order $6$. Then $\gamma(\Gamma_G)=2$ if and only if
the intersection of
any two cyclic subgroups of order $6$ is of order $3$.
\end{lem}
\bpf
Suppose that $\gamma(\Gamma_G)=2$. If the intersection of any two cyclic subgroups of order $6$ is of order $1$ or $2$, then it is easy to see that $\Gamma_G$ has a subgraph isomorphic to $K_1+3K_4$ that has genus three by Theorem \ref{moo}, and so a contradiction. Therefore, there exist at least two cyclic subgroups of order $6$ in $G$ such that their intersection has order $3$. Now
by Lemma \ref{claim}, we get the desired result.

For the converse, let $H_1,H_2$ and $H_3$ be the three cyclic subgroups of $G$.
We assume that $H_1=\{e,f_1,f_2,x,g_1,g_2\}$, $H_2=\{e,f_3,f_4,y,g_1,g_2\}$ and $H_3=\{e,f_5,f_6,z,g_1,g_2\}$, where $|f_i|=6,|g_j|=3$ and $|x|=|y|=|z|=2$ for $i=1,\ldots,6$ and $j=1,2$.

Let $\Gamma$ be the subgraph of $\Gamma_G$ induced by $\cup_{i=1}^3H_i$.
Denote by $\Gamma'$ the subgraph of $\Gamma$ obtained by deleting the vertices $x,y,z$ and
the edges $eg_1,eg_2,g_1g_2,f_1f_2,f_3f_4,f_5f_6$ in $\Gamma$.
Then $\Gamma'$ is isomorphic to $K_{3,6}$.
Note that $\gamma(\Gamma')=1$, the boundary of each face is a $4$-cycle when drawing $\Gamma'$ without crossings on a torus, and any two faces of $\Gamma'$ have at most one boundary edge in common (see \cite[Remark 1.4]{CSW10}).
Now we proceed to prove $\gamma(\Gamma)\ge 2$ by a deletion and insertion argument.

Suppose, to the contrary, that $\gamma(\Gamma)=1$.
Then $\Gamma$ can be embedded in a torus without crossings.
Fix an embedding $\mathbb{E}$ of $\Gamma$ on a tours. By deleting some vertices and edges from the embedding of $\Gamma$, we can get an embedding $\mathbb{E}'$ of $\Gamma'$ on a tours. This implies that
all faces of $\mathbb{E}$ can be recovered by inserting some vertices and edges
into the faces of $\mathbb{E}'$.
Let $F'$ be the face of $\mathbb{E}'$ into which $x$ is inserted during the recovering process from $\mathbb{E}'$ to $\mathbb{E}$. Note that $xf_1,xf_2\in E(\Gamma)$.
We obtain the face $F'$ as shown in Figure \ref{f2}, where $\{u,v\}\subseteq \{e,g_1,g_2\}$.
\begin{figure}[hptb]
  \centering
  \includegraphics[width=6cm]{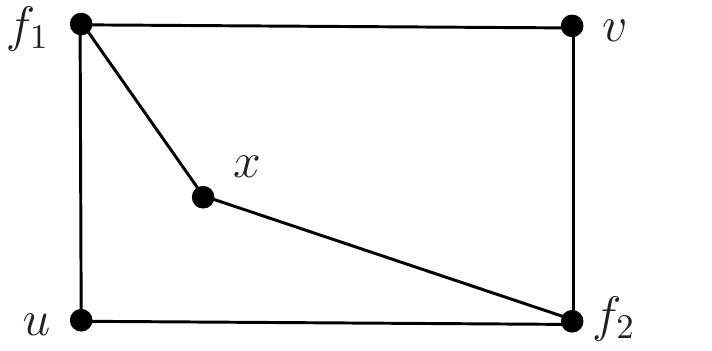}\\
  \caption{Insert $a$ in $F'$.}\label{f2}
\end{figure}
Since $\{e,g_1,g_2\}$ is a clique, one has that $u$ is adjacent to $v$ in $\Gamma$.
Thus, in order to insert the edge $uv$  without crossings, $\mathbb{E}'$ must have a face  different with $F'$ such that its boundary is a $4$-cycle containing $u$ and $v$. This implies that there are two faces of $\mathbb{E}'$ such that they have two boundary edges in common, a contradiction. Thus, $\gamma(\Gamma)\ge 2$.

On the other hand, we can embed $\Gamma$ into $\mathbb{S}_2$ as shown in Figure \ref{fn2}. This means that $\gamma(\Gamma)=2$. Now note that $\pi_e(G\setminus \cup_{i=1}^3H_i)\subseteq \{2,3,4\}$. By Lemma \ref{condi} one has that $\gamma(\Gamma_G)=2$.
\begin{figure}[hptb]
  \centering
  \includegraphics[width=7cm]{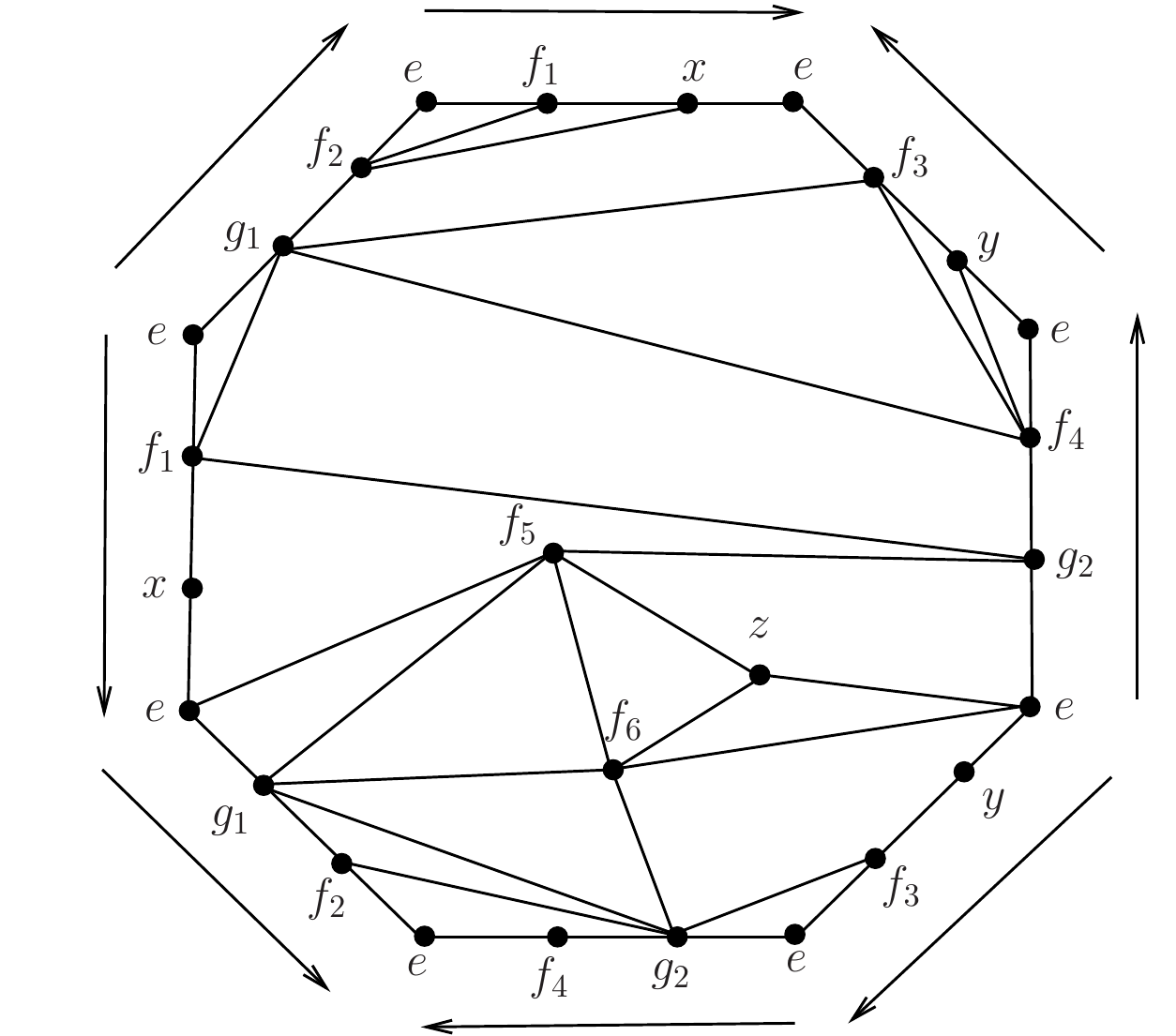}\\
  \caption{An embedding of $\Gamma$ on $\mathbb{S}_2$.}\label{fn2}
\end{figure}
\epf

\begin{lem}\label{four6}
Let $G$ be a group with $\pi_e(G)\subseteq\{1,2,3,4,6\}$. Suppose that $G$ has precisely four cyclic subgroups $H_1,\ldots,H_4$ of order $6$. If $\gamma(\Gamma_G)=2$, then
$|H_i\cap H_j|=3$ and $|H_s\cap H_t|=3$, where $\{i,j,s,t\}=\{1,2,3,4\}$.
\end{lem}
\bpf
By considering all possible cases and  Theorem~\ref{moo}, it is not difficult to get the desired result.
\epf

\begin{lem}\label{five6}
Let $G$ be a group with $\pi_e(G)\subseteq\{1,2,3,4,6\}$. If $G$ has at least five cyclic subgroups $H_1,\ldots,H_5$ of order $6$, then $\gamma(\Gamma_G)\ge 3$.
\end{lem}
\bpf
If $|\cap_{i=1}^5H_i|\ne 3$, we check all possible cases and then by Theorem~\ref{moo},
it is easy to see that $\gamma(\Gamma_G)\ge 3$. Thus, we may suppose that $|\cap_{i=1}^5H_i|= 3$. Let $H_1=\{e,f_1,f_2,a,g_1,g_2\}$, $H_2=\{e,f_3,f_4,b,g_1,g_2\}$, $H_3=\{e,f_5,f_6,c,g_1,g_2\}$, $H_4=\{e,f_7,f_8,d,g_1,g_2\}$ and $H_5=\{e,f_9,f_0,w,g_1,g_2\}$,
where $|a|=|b|=|c|=|d|=|w|=2$, $|f_i|=6$ and $|g_j|=3$ for $i=0,\ldots,9$ and $j=1,2$.
Denote by $\Delta$ the subgraph of $\Gamma_G$ induced by the vertices $\{e,f_i,g_1,g_2: i=0,\ldots,9\}$. Write $\Delta'=\Delta-\{f_1f_2,f_3f_4,f_5f_6,f_7f_8,eg_1,eg_2,g_1g_2\}$ and
$\Delta''=\Delta'-\{f_9,f_0\}$. Then $\Delta''\cong K_{3,8}$. Since $\gamma(K_{3,8})=2$, one has
$\gamma(\Delta)\ge 2$.
Next we prove $\gamma(\Delta)\ge 3$ by a deletion and insertion argument.

Assume, for a contradiction, that $\gamma(\Delta)= 2$. Then $\gamma(\Delta')= 2$.
Fix an  embedding $\mathbb{E}'$ of $\Delta'$ on $\mathbb{S}_2$.
Thus, we may get an embedding $\mathbb{E}''$ of $\Delta''$ on $\mathbb{S}_2$ such that
$\mathbb{E}'$ can be recovered by inserting $f_9,f_0$ and all edges
incident with $f_9,f_0$ into the embedding $\mathbb{E}''$.
Since $\gamma(\Delta'')=2$, by Theorem \ref{ef} there are $11$ faces when  drawing $\Delta''$ without crossings on $\mathbb{S}_2$ and thereby, each of the faces is  a $4$-cycle or a $6$-cycle.
Let $F_{t}''$ denote the face of $\Delta''$ into which $f_0$ is inserted during the recovering process from $\mathbb{E}''$ to $\mathbb{E}'$. Since $f_9$ is adjacent to $f_0$ in $\Delta'$, $f_9$ should be also inserted into $F_{t}''$ to avoid any crossing. Moreover, since $ef_i,g_1f_i,g_2f_i\in E(\Delta')$ for $i=0,9$, one has that $e,g_1,g_2$ lie in the boundary of $F_{t}''$. This implies that $F_{t}''$ is a $6$-cycle. Then after
inserting $f_0,f_9$ and $f_0f_9,ef_0,ef_9,g_1f_0,g_1f_9,f_9g_2$ into $F_{t}''$ we get Figure~\ref{f4} as below.
\begin{figure}[hptb]
  \centering
  \includegraphics[width=7cm]{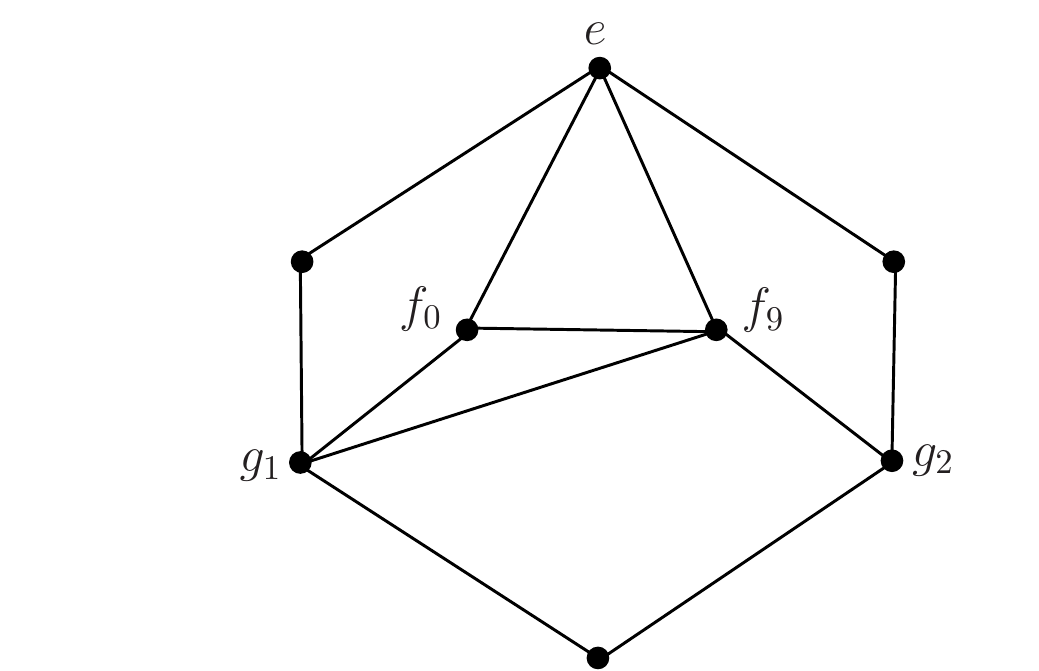}\\
  \caption{The face $F_{t}''$.}\label{f4}
\end{figure}
But, it is easy to see from
Figure~\ref{f4} that we can not insert the edge $f_0g_2$ into $F_{t}''$ without crossings, a contradiction. Thus, we conclude that $\gamma(\Delta)\ne 2$ and so $\gamma(\Delta)\ge 3$.
It follows that $\gamma(\Gamma_G)\ge 3$.
\epf

{\noindent \em Proof of Theorem \ref{mainthm1}.}
Suppose that $\gamma(\Gamma_{G})=2$. Then by Lemma \ref{cycg} we see that
$$\pi_e(G)\subseteq \{1,2,3,4,5,6,7,8\}.$$

Suppose that $P\ne 1$ is a Sylow $5$-subgroup of $G$. Then $\gamma(\Gamma_P)\le \gamma(\Gamma_{G})=2$.
If $|P|\ge 25$, then $G$ has a subgroup $K$ isomorphic to $\mathbb{Z}_5\times \mathbb{Z}_5$
and so $\gamma(\Gamma_K)=6$ by
Theorem \ref{moo}, a contradiction.
It follows that $P\cong \mathbb{Z}_5$.
Furthermore, by
Sylow's Theorem, the number of Sylow $5$-subgroups of $G$ is $5k+1$ for some integer $k$.
Note that by Theorem \ref{moo}, the subgraph of $\Gamma_{G}$ induced by all Sylow $5$-subgroups of $G$ has genus $5k+1$. Therefore, we have that $k=0$ and so $P$ is normal in $G$.
Since $\gamma(\Gamma_{\mathbb{Z}_8})=2$, $G$ has no elements of order $8$ by Theorem \ref{moo}.
Note that every group of order $30$ has an element of order $15$, and every group of order $35$ is isomorphic to $\mathbb{Z}_{35}$. So $G$ has no elements of order $7$ and $6$.
Consequently, we conclude that $\pi_e(G)\subseteq \{1,2,3,4,5\}$. Since $G$ has precisely one subgroup of order $5$ and $\gamma(\Gamma_P)=1$, one has $\gamma(\Gamma_{G})=1$ by Lemma \ref{condi}, a contradiction.
This implies that $5\notin \pi(G)$.
Similarly, we have that $7\notin \pi(G)$.

Suppose that $G$ has an element $g$ of order $8$. Then it follows from Theorem \ref{moo} that $6\notin \pi_e(G)$ and $G$ has only one cyclic subgroup of order $8$. If there exists an element
$a$ in $G$ such that $|a|=3$, then $\langle a,g\rangle$ is of order $24$ and $\pi_e(\langle a,g\rangle)=\{1,2,3,4,8\}$, a contradiction since such a group does not exist.
Therefore, in this case $G$ is a $2$-group. By Lemma \ref{p2g} we get the required result.

Now we may assume that $\pi_e(G)\subseteq\{1,2,3,4,6\}$. By Theorem \ref{plan} $G$ has some elements of order $6$.
If $G$ has exactly one cyclic subgroup of order $6$, then
by Lemma \ref{condi} it is easy to see that $\Gamma_G$ has genus one. Therefore, we conclude that
$G$ has at least two cyclic subgroups of order $6$.
Now  the necessity follows from Theorems \ref{gtheory3} and \ref{gtheory4}, and Lemmas \ref{two6}, \ref{three6}, \ref{four6} and \ref{five6}.

For the converse, it follows from  Theorems \ref{gtheory3}, Lemmas \ref{p2g} and \ref{three6}. \qed

\section{Proof of Theorem \ref{mainthm2}}

In this section we show Theorem \ref{mainthm2}. We begin with the following lemma.

\begin{lem}\label{ng21}
$\overline{\gamma}(\Gamma_{\mathbb{Z}_n})\ne 2$ for each $n$. In particular,
$\overline{\gamma}(\Gamma_{\mathbb{Z}_n})\ge 3$ for $n\ge 7$.
\end{lem}
\bpf
Note that if $n\ge 7$, then the clique number of $\Gamma_{\mathbb{Z}_n}$ is greater than or equal to $7$. Now the result follows from Theorem \ref{ccgenus}.
\epf

\begin{lem}\label{ng22}
Let $G$ be a group with $\pi_e(G)\subseteq\{1,2,3,4,6\}$. Suppose that $G$ has at least three cyclic subgroups of order $6$. Then $\overline{\gamma}(\Gamma_G)> 2$.
\end{lem}
\bpf
Let $H_1,H_2,H_3$ be three cyclic subgroups of order $6$ of $G$.
If $|H_s\cap H_t|\ne 3$ for each two $s,t$, then $\Gamma_G$ has a subgraph isomorphic to
$K_1+(K_4\cup K_4 \cup K_4)$ that has nonorientable genus three by Theorem \ref{moo} and hence $\overline{\gamma}(\Gamma_G)> 2$, as desired. Thus, we now may assume that there exist two cyclic subgroups of order $6$ in $G$ such that their intersection is of order $3$. Without loss of generality, let $|H_1\cap H_2|=3$.
Denote by $\Delta$ the subgraph induced by $H_1\cup H_2$. Next we prove that $\overline{\gamma}(\Delta)=2$.

Let $H_1=\{e,g_1,g_2,f_1,f_2,a\}$ and $H_1=\{e,g_1,g_2,f_3,f_4,b\}$, where $|g_1|=|g_2|=3$, $|a|=|b|=2$ and $|f_i|=6$ for $i=1,\ldots, 4$. Then $\Delta- \{a,b\}$ is the graph as shown in Figure~\ref{a1}, which is isomorphic to the graph $B_1$ in \cite{GHW}.
By the main result of \cite{GHW}, one has $\overline{\gamma}(\Delta)\ge 2$.
\begin{figure}[hptb]
  \centering
  \includegraphics[width=7cm]{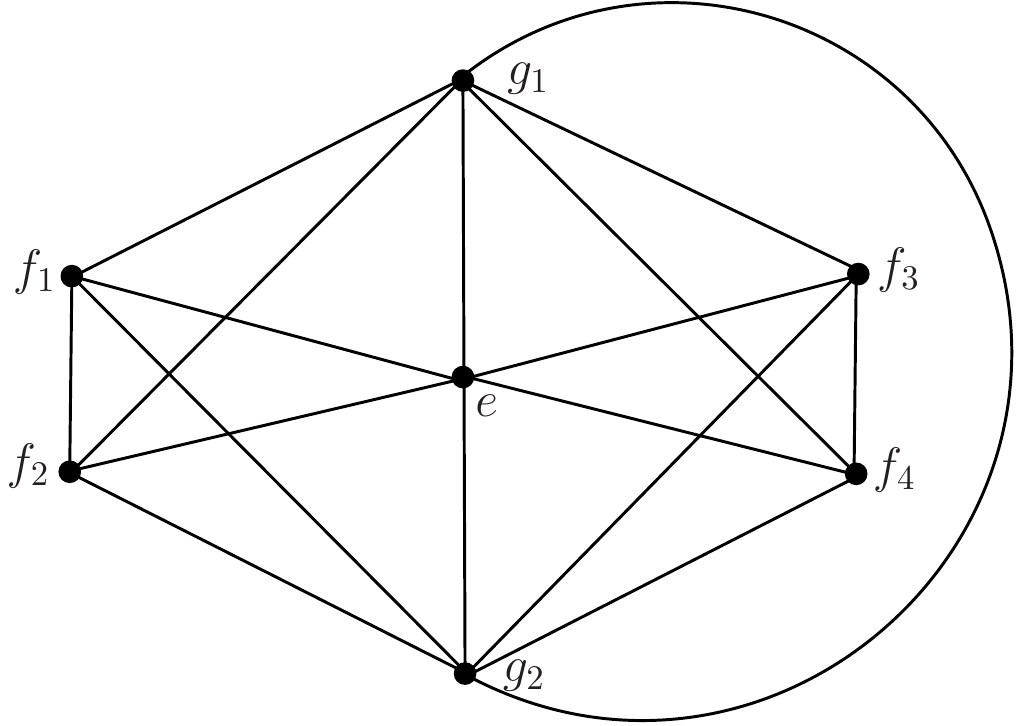}\\
  \caption{An obstruction for the projective plane.}\label{a1}
\end{figure}
On the other hand, we can embed $\Delta$ on
$\mathbb{N}_2$ as shown in Figure~\ref{a2}, so $\overline{\gamma}(\Delta)=2$.
\begin{figure}[hptb]
  \centering
  \includegraphics[width=7cm]{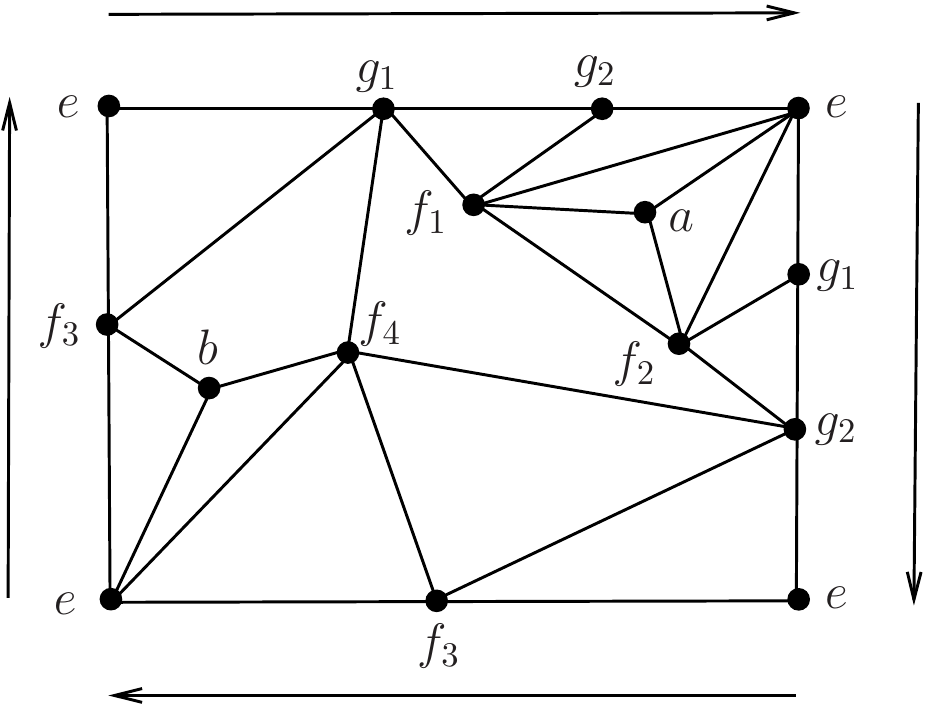}\\
  \caption{An embedding of $\Delta$ on $\mathbb{N}_2$}\label{a2}
\end{figure}

Suppose that $|H_3\cap H_1|\ne 3$. Then $\Gamma_G$ contains a subgraph $\Gamma$ that has two
blocks $\Delta$ and a graph isomorphic to $K_5$. Clearly, $\gamma(\Delta)\ge 1$.
Since we can embed $\Delta$ on a torus as shown in Figure~\ref{f3}, one has $\gamma(\Delta)=1$.
Now by Theorem \ref{moo} we have that $\overline{\gamma}(\Gamma)=3$ and thereby
$\overline{\gamma}(\Gamma_G)\ge 3$, as required.
\begin{figure}[hptb]
  \centering
  \includegraphics[width=7cm]{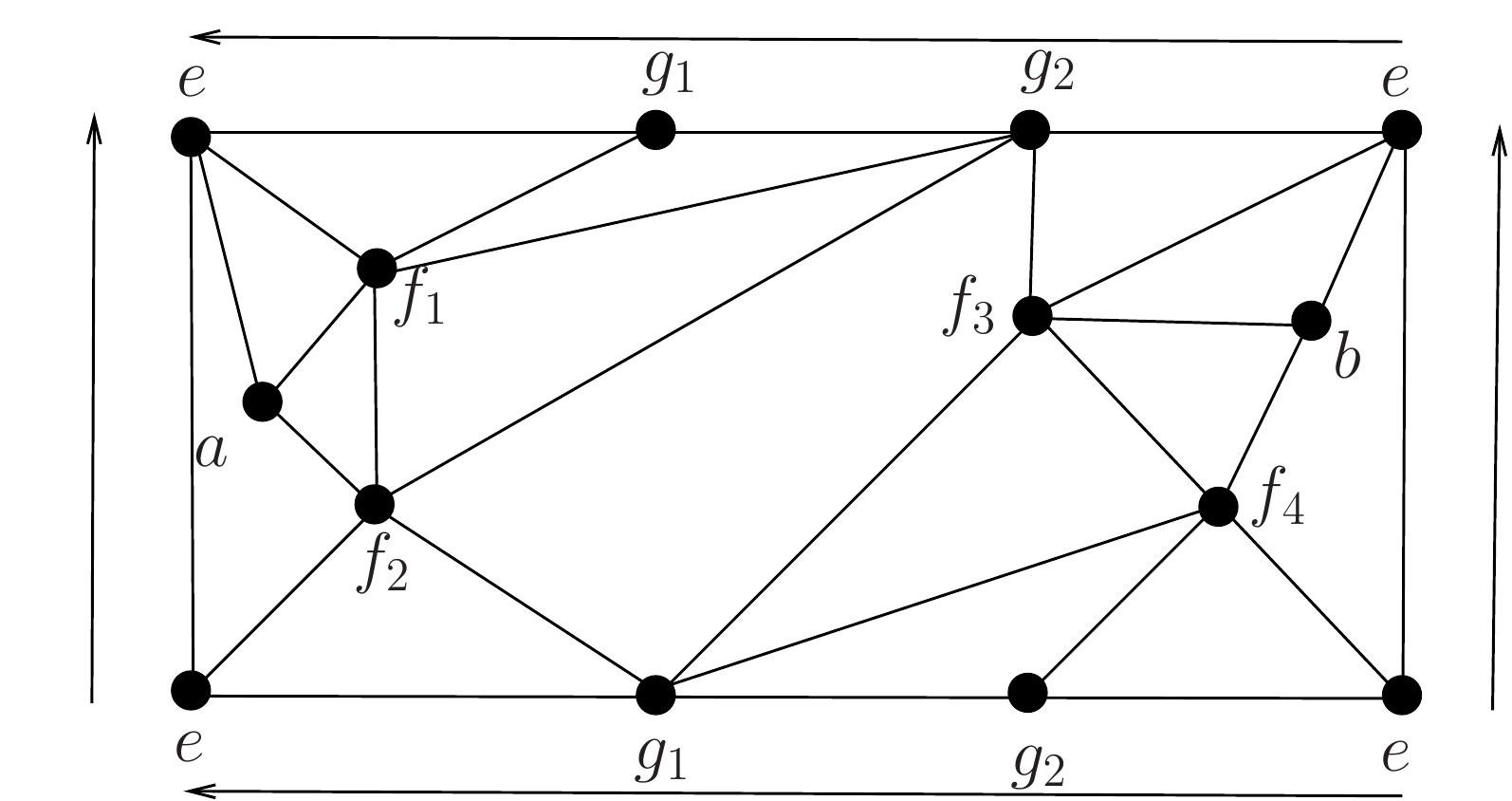}\\
  \caption{An embedding of $\Delta$ on a torus}\label{f3}
\end{figure}
Therefore, we may suppose that $|H_3\cap H_1|=3$. Then $|H_1\cap H_2 \cap H_3|=3$. Write
$H_3=\{e,g_1,g_2,f_5,f_6,c\}$ where $|f_5|=|f_6|=6$ and $|c|=2$. Let $\Lambda$ denote the subgraph induced by $\cup_{i=1}^3 H_i$. Then $\Lambda-\{f_5,f_6,c\}=\Delta$.
It is clear that $\overline{\gamma}(\Lambda)\ge \overline{\gamma}(\Delta)= 2$. In order to end the proof, next we show $\overline{\gamma}(\Lambda)\ge 3$.

Suppose, towards a contradiction, that $\overline{\gamma}(\Lambda)=2$.
By Theorem \ref{ef} there are $15$ faces when drawing $\Lambda$ without crossings on $\mathbb{N}_2$. Fix an  embedding
of $\Lambda$ on $\mathbb{N}_2$ and let $\{F_1,\ldots,F_{15}\}$ be the set of all faces of
$\Lambda$ corresponding to this embedding.
So the faces obtained by deleting $f_5,f_6,c$ and all edges incident with them from $\{F_1,\ldots,F_{15}\}$ form an embedding of $\Delta$ on $\mathbb{N}_2$. In other words,
the faces $\{F_1,\ldots,F_{15}\}$ can be recovered by inserting $f_5,f_6,c$ and all edges
incident with them into an embedding of $\Delta$ on $\mathbb{N}_2$.
We first insert $f_5$ into a face $F$. Since $f_5$ and $f_6$  are adjacent in $\Lambda$, $f_6$ must be inserted into $F$. Note that $\{f_5e,f_5g_1,f_5g_2\}\subseteq E(\Lambda)$.
So $e,g_1,g_2$ lie in the boundary of $F$.  However, we can not  insert the edges $f_6e,f_6g_1,f_6g_2$ into $F$ without crossings (similarly, cf. Figure \ref{f4}), a contradiction.
\epf

\bigskip

{\noindent \em Proof of Theorem \ref{mainthm2}.}
Suppose for a contradiction that $\overline{\gamma}(\Gamma_{G})=2$ for some group $G$. Then by Lemma \ref{ng21} we have $\pi_e(G)\subseteq \{1,2,3,4,5,6\}$.

Suppose that $G$ has a Sylow $5$-subgroup $P$.
If $|P|\ge 25$, then $\Gamma_G$ has a subgraph isomorphic to $\Gamma_{\mathbb{Z}_5\times \mathbb{Z}_5}$ that has $6$ blocks and each of its blocks is isomorphic to $K_5$, and so $\overline{\gamma}(\Gamma_{\mathbb{Z}_5\times \mathbb{Z}_5})=6$ by Theorem \ref{moo}, a contradiction.
Therefore, we have $|P|=5$. Considering the subgraph induced by all Sylow $5$-subgroups, similarly, we have that $G$ has a unique Sylow $5$-subgroup. This implies that $P$ is normal in $G$.
If $G$ has an element of order $6$, then $G$ has a subgroup of order $30$, a contradiction since
every group of order $30$ has an element of order $15$. Thus, in this case one has
$\pi_e(G)\subseteq \{1,2,3,4,5\}$. Note that $\pi_e(G\setminus P)\subseteq \{2,3,4\}$.
In view of Lemma \ref{condi}, one has  $\overline{\gamma}(\Gamma_{G})=1$, a contradiction.

Thus, now we may assume that $\pi_e(G)\subseteq \{1,2,3,4,6\}$. If $G$ has precisely one cyclic
subgroup of order $6$, then by Lemma \ref{condi} it is easy to see that $\overline{\gamma}(\Gamma_{G})=1$, a contradiction.
So $G$ has at least two cyclic subgroups of order $6$.
Now we can get the final contradiction by
Lemmas~\ref{two6} and \ref{ng22}.  \qed

\section*{Acknowledgement}
K. Wang's research is supported by National Natural Science Foundation of China (11271047,
11371204) and the Fundamental Research Funds for the Central University of China.

\end{document}